\documentclass[11pt]{article}
\usepackage{latexsym, amssymb}

\def \F {{\mathbb F}}

\def \Tr {{\rm Tr_n}}
\newtheorem{theorem}{Theorem}

\newtheorem{lemma}{Lemma}

\newtheorem{remark}{Remark}
\newtheorem{corollary}{Corollary}
\def\bproof{\textbf{Proof}: }
\def\eproof{\hfill$\Box$}
\begin{document}

\begin{center}
{\Large \bf A Construction of Weakly and Non-Weakly Regular Bent Functions}
\end{center}
\begin{center}
Ay\c ca \c Ce\c smelio\u glu$^{*}$, Gary McGuire$^{\dagger}$,\footnote{Research supported by the Claude
Shannon Institute, Science Foundation Ireland Grant 06/MI/006} Wilfried Meidl$^{*}$
\end{center}
$^{*}$ Sabanc{\i} University, MDBF, Orhanl\i, {\rm 34956} Tuzla, \.Istanbul, Turkey. \\
$^{\dagger}$ School of Mathematical Sciences, University College Dublin, Ireland

\begin{abstract}
In this article a technique for constructing $p$-ary bent functions
from near-bent functions is presented. Two classes of quadratic
$p$-ary functions are shown to be near-bent. Applying the 
construction of bent functions to these classes of near-bent functions
yields classes of non-quadratic bent functions. We show that one construction
in even dimension yields weakly regular bent functions. For other constructions,
we obtain both weakly regular and non-weakly regular bent functions. In particular
we present the first known infinite class of non-weakly regular bent functions.

{\it Keywords:} Bent function, Near-bent function, Semi-bent function, Weakly regular, Non-weakly regular, Fourier transform.
\end{abstract}

\section{Introduction}


Let $p$ be a prime, and let $V_n$ be any $n$-dimensional vector space over $\F_p$.
For a function $f$ from $V_n$ to $\F_p$ the {\it Fourier transform (or Walsh transform)} of
$f$ is the complex valued function $\widehat{f}$ on $V_n$ given by
\[ \widehat{f}(b) = \sum_{x\in V_n}\epsilon_p^{f(x)-\langle b,x\rangle} \]
where $\epsilon_p = e^{2\pi i/p}$ and $\langle,\rangle$ denotes any inner 
product on $V_n$. 
We say that $f$ is a
{\it  bent} function if $|\widehat{f}(b)|^2 = p^n$ for all $b
\in V_n$.  If $p=2$ then $\epsilon_p=-1$ and $\widehat{f}(b)$ is an integer,
so a necessary condition for the existence of a bent function is that $n$ is even.
This does not hold for odd $p$, where bent functions can exist for both odd and even $n$.
When $p$ is odd, bent functions are sometimes called $p$-ary bent functions.

As all vector spaces of dimension $n$ over $\F_p$ are isomorphic we may 
associate $V_n$ with the finite field $\F_{p^n}$. We then usually use the
inner product $\langle x,y\rangle = \Tr(xy)$ where $\Tr(z)$ denotes the 
absolute trace of $z\in\F_{p^n}$. In this framework the Fourier transform 
of a function $f$ from $\F_{p^n}$ to $\F_p$ is the complex valued function 
on $\F_{p^n}$ given by
\[ \widehat{f}(b) = \sum_{x\in\F_{p^n}}\epsilon_p^{f(x)-\Tr(bx)}. \]

Often one considers the 
{\it normalized} Fourier coefficient
$p^{-n/2}\widehat{f}(b)$ of a  bent function.
For any $p$, we can only say a priori that the normalized Fourier coefficients lie on
the unit circle. 
For $p=2$, a bent function must have normalized Fourier coefficients
$\pm 1$, because the Fourier coefficients  are real.
For odd $p$, the values of the normalized Fourier coefficients
of a bent function are also quite constrained,
(cf. \cite{hk}, \cite[Property 8]{ksw}).
The possibilities are as follows:
\begin{equation}
\label{(2)} p^{-n/2}\widehat{f}(b) =
\left\{\begin{array}{r@{\quad:\quad}l}
\pm \epsilon_p^{f^*(b)} & n\;\mbox{even or}\;n\;\mbox{odd and}\;p \equiv 1\bmod 4 \\
\pm i\epsilon_p^{f^*(b)} & n\;\mbox{odd and}\;p \equiv 3\bmod 4
\end{array}\right.
\end{equation}
where $f^*$ is a function from $\F_{p^n}$ to $\F_p$ that
by definition gives the exponent of $\epsilon_p$. 

A bent function $f$ is called {\it regular} if, for all $b \in
\F_{p^n}$, we have 
\[
p^{-n/2}\widehat{f}(b) =\epsilon_p^{f^*(b)},
\]
i.e., the coefficient of $\epsilon_p^{f^*(b)}$ is always $+1$.
Thus the normalized Fourier coefficients of a regular bent function are
a subset of (in fact the full set of)
the $p$-th roots of unity. 
This is a natural generalization of the binary situation, where the normalized
Fourier coefficients are $\pm 1$.
 It is obvious from equation $(\ref{(2)})$
that regular bent functions can only exist for even $n$ and for odd $n$
with $p\equiv 1\bmod 4$.
For example, when $p=3$,  regular bent functions can only exist
in even dimensions, and
the normalized Fourier coefficients are shown in 
the complex plane in Fig 1.

\setlength{\unitlength}{0.45mm}

$$\begin{picture}(80,90)
\put(40,83){\line(0,-1){80}}
\put(0,40){\line(1,0){80}}
\put(80,40){\circle*{3}}
\put(20,70){\circle*{3}}
\put(20,10){\circle*{3}}
\put(80,32) {$\scriptstyle{1}$}
\put(20,62) {$\scriptstyle{\epsilon_3}$}
\put(20,15) {$\scriptstyle{\epsilon_3^2}$}
\put(40,-5){$\textrm{Fig\ 1, $p=3$, even dimension}$}
\end{picture}
$$

\bigskip

A bent function $f$ is called {\it weakly regular} if, for all $b \in
\F_{p^n}$, we have 
\[
p^{-n/2}\widehat{f}(b) =\zeta\ \epsilon_p^{f^*(b)}
\]
for some complex number $\zeta$ with absolute value $1$   
(see \cite{ksw}). 
By (\ref{(2)}), $\zeta$ can only be $\pm 1$ or $\pm i$.
Thus the normalized Fourier coefficients of a weakly regular bent function are
a rotation (through some multiple of $\pi/2$) 
of the $p$-th roots of unity.
For $p=3$, the three possibilities for the normalized Fourier coefficients are shown in 
Fig 2, Fig 3 and Fig 4.
In this ternary case, weakly regular bent functions 
(that are not regular)
can  exist
in both odd and even dimensions.

$$\begin{picture}(80,90)
\put(40,83){\line(0,-1){80}}
\put(0,40){\line(1,0){80}}
\put(0,40){\circle*{3}}
\put(60,70){\circle*{3}}
\put(60,10){\circle*{3}}
\put(0,32) {$\scriptstyle{-1}$}
\put(60,62) {$\scriptstyle{-\epsilon_3^2}$}
\put(60,15) {$\scriptstyle{-\epsilon_3}$}
\put(40,-5){$\textrm{Fig\ 2,  $p=3$, even dimension}$}
\end{picture}
$$
\bigskip

$$\begin{picture}(80,90)
\put(40,83){\line(0,-1){80}}
\put(0,40){\line(1,0){80}}
\put(10,20){\circle*{3}}
\put(40,79){\circle*{3}}
\put(70,20){\circle*{3}}
\put(43,80) {$\scriptstyle{i}$}
\put(10,23) {$\scriptstyle{i \epsilon_3}$}
\put(70,23) {$\scriptstyle{i\epsilon_3^2}$}
\put(40,-5){$\textrm{Fig\ 3,  $p=3$, odd dimension}$}
\end{picture}
$$
\bigskip

$$\begin{picture}(80,90)
\put(40,83){\line(0,-1){80}}
\put(0,40){\line(1,0){80}}
\put(10,60){\circle*{3}}
\put(40,4){\circle*{3}}
\put(70,60){\circle*{3}}
\put(43,5) {$\scriptstyle{-i}$}
\put(10,63) {$\scriptstyle{-i \epsilon_3^2}$}
\put(70,63) {$\scriptstyle{-i\epsilon_3}$}
\put(40,-5){$\textrm{Fig\ 4,  $p=3$, odd dimension}$}
\end{picture}
$$
\bigskip

Almost all known $p$-ary bent functions are weakly regular.
Until this paper, there are just a few sporadic examples of non-weakly regular
bent functions known (see \cite{hk,hk1}).
If $p=3$, a non-weakly regular bent function in even dimensions would have
normalized Fourier coefficients as in Fig 5
(Fig 1 and Fig 2 combined).

\bigskip
  
$$\begin{picture}(80,90)
\put(40,83){\line(0,-1){80}}
\put(0,40){\line(1,0){80}}
\put(0,40){\circle*{3}}
\put(60,70){\circle*{3}}
\put(60,10){\circle*{3}}
\put(80,40){\circle*{3}}
\put(20,70){\circle*{3}}
\put(20,10){\circle*{3}}
\put(40,-5){$\textrm{Fig\ 5, $p=3$,  even dimension}$}
\end{picture}
$$
\bigskip

A ternary non-weakly regular bent function in odd dimensions would have
normalized Fourier coefficients as in Fig 6
(Fig 3 and Fig 4 combined).
We give a construction of such functions in this paper (for any $p$).

$$\begin{picture}(80,90)
\put(40,83){\line(0,-1){80}}
\put(0,40){\line(1,0){80}}
\put(10,20){\circle*{3}}
\put(40,79){\circle*{3}}
\put(70,20){\circle*{3}}
\put(10,60){\circle*{3}}
\put(40,4){\circle*{3}}
\put(70,60){\circle*{3}}
\put(40,-5){$\textrm{Fig\ 6, $p=3$,  odd dimension}$}
\end{picture}
$$
\bigskip

For the binary case, where bent functions in odd dimension do not
exist, the notion of {\it near-bent} functions was introduced in
\cite{lg}. We generalize this now to characteristic $p$,
and we call a function $f$ from $\F_{p^n}$ to $\F_p$
\emph{near-bent} if, for all $b \in \F_{p^n}$,  $|\widehat{f}(b)|^2 =
p^{n+1}$ or $0$. We remark that in \cite{cpt} the
term semi-bent function is used for the same concept in characteristic 2.

In this article we first generalize to characteristic $p$  
the technique presented in \cite{lg}
(see also \cite{cpt}) that constructed binary bent functions from
near-bent functions. In Section \ref{sec2} we
illustrate the principle of the construction. 
The idea is to choose near-bent functions in dimension $n$,
which will be mappings from $\F_{p^n}$ to $\F_p$,
and glue them together using another copy of $\F_p$ to obtain
bent functions in dimension $n+1$. 
These functions will be from
$\F_{p^n}\times \F_p$ to $\F_p$.
The near-bent functions must be chosen so that
the supports of their Fourier transforms are disjoint,
so exactly one of their Fourier transforms is nonzero at any point.
We will show that varying the glueing coefficients from the extra copy of $\F_p$
can drastically change the nature of the bent function,
thereby demonstrating that inequivalent bent functions can be
constructed by this method. For example,
both weakly regular and non-weakly regular bent functions
can be found by a simple tweak of the $\F_p$-coefficients 
(which are denoted $c_k$ below,
see Theorem \ref{sum}).

In Section \ref{sec3}
we present a class of $p$-ary quadratic binomials that are near-bent.
Another feature of our construction is that we obtain infinite families of 
\emph{non}-quadratic bent
functions from quadratic near-bent functions, so we are able to leave
the quadratic world. 
The construction of bent
functions from quadratic near-bent functions is described in detail
in Section \ref{sec4}. 
In Section \ref{sec5} we prove that one class of bent functions
obtained with the near-bent functions introduced in Section
\ref{sec3} is weakly regular. For the general construction
we give if and only if conditions to yield weakly regular bent functions.
Using these conditions we are then able to give infinite classes of 
non-weakly regular bent functions.

\section{Obtaining bent functions from near-bent functions}
\label{sec2}

Let $f$ be a function from $\F_{p^n}$ to $\F_p$, and $\widehat{f}$
denote its Fourier transform. The support of $\widehat{f}$ is then
defined to be the set $supp(\widehat{f}) =
\{b\in\F_{p^n}\;|\;\widehat{f}(b)\ne 0\}$.
For any $p$-ary function $f$ we have
\begin{eqnarray*}
\sum_{b\in\F_{p^n}}\left|\widehat{f}(b)\right|^2&=&\sum_{b\in\F_{p^n}}\sum_{x,y \in\F_{p^n}}\epsilon_p^{f(x)-\Tr(bx)-(f(y)-\Tr(by))}\\
&=&\sum_{x,y\in\F_{p^n}}\epsilon_p^{f(x)-f(y)}\sum_{b
\in\F_{p^n}}\epsilon_p^{\Tr(b(y-x))}.
\end{eqnarray*}
Observing that $\sum_{b \in\F_{p^n}}\epsilon_p^{\Tr(b(y-x))}=0$ if
$x\neq y$, and $\sum_{b\in\F_{p^n}}\epsilon_p^{\Tr(b(y-x))}=p^n$ if
$x=y$, we obtain the special case of {\it Parseval's relation}:
\[\sum_{b\in\F_{p^n}}\left|\widehat{f}(b)\right|^2=\sum_{x,y\in\F_{p^n},x=y}p^n=p^{2n}.\]
For a near-bent function $f$, clearly
\[\sum_{b\in\F_{p^n}}\left|\widehat{f}(b)\right|^2=\left|supp(\widehat{f})\right|p^{n+1}\]
and combining this with Parseval's relation gives
\[\left|supp(\widehat{f})\right|=p^{n-1}.\]

The following theorem presents how to obtain $p$-ary bent functions
from a set of $p$ near-bent functions
$f_0(x),f_1(x),\cdots,f_{p-1}(x)$ from $\F_{p^n}$ to $\F_p$ with
$supp(\widehat{f_i})\cap supp(\widehat{f_j})=\emptyset$ for $i\neq
j$. We remark that then
$\bigcup_{i=0}^{p-1}supp(\widehat{f_i})=\F_{p^n}$.
\begin{theorem}
\label{thm1} Let $f_0(x),f_1(x),\cdots,f_{p-1}(x)$ be near-bent
functions from $\F_{p^n}$ to $\F_p$ such that
$supp(\widehat{f_i})\cap supp(\widehat{f_j})=\emptyset$ for $0 \le
i\neq j \le p-1$. Then the function $F(x,y)$ from
$\F_{p^n}\times\F_p$ to $\F_p$ defined by
\[ F(x,y)=(p-1)\sum_{k=0}^{p-1}\frac{y(y-1)\cdots(y-(p-1))}{y-k}f_k(x) \]
is bent.
\end{theorem}

\bproof For $(a,b), (x,y) \in \F_{p^n}\times\F_p$ 
the inner product we use is $\Tr(ax)+by$.
The Fourier transform
$\widehat{F}$ of $F$ at $(a,b)$ is
\begin{eqnarray*}
\widehat{F}(a,b)&=&\sum_{x \in \F_{p^n}, y \in \F_p}\epsilon_p^{F(x,y)-\Tr(ax)-by}\\
&=&\sum_{y \in \F_p}\epsilon_p^{-by}\sum_{x \in \F_{p^n}}\epsilon_p^{F(x,y)-\Tr(ax)}\\
&=&\sum_{y \in \F_p}\epsilon_p^{-by}\sum_{x \in \F_p^n}\epsilon_p^{(p-1)!(p-1)f_y(x)-\Tr(ax)}\\
&=&\sum_{y \in \F_p}\epsilon_p^{-by}\sum_{x \in \F_p^n}\epsilon_p^{f_y(x)-\Tr(ax)}\\
&=&\sum_{y \in \F_p}\epsilon_p^{-by}\widehat{f_y}(a).
\end{eqnarray*}
As each $a \in \F_{p^n}$ belongs to the support of exactly one
$\widehat{f_y}$, $y \in \F_p$, for this $y$ we have
$\left|\widehat{F}(a,b)\right|=|\epsilon_p^{-by}\widehat{f_y}(a)|=p^{\frac{n+1}{2}}$. \eproof

\section{Monomial and binomial quadratic near-bent functions}
\label{sec3}

Recall that a function $f$ from $\F_{p^n}$ to $\F_p$ of the form
\[ f(x) = \Tr\left(\sum_{i=0}^la_ix^{p^i+1}\right) \]
is called {\it quadratic}, its algebraic degree is two (if $f$ is
not constant), see \cite{cpt,hk}. The following theorem giving
the Fourier spectrum of a quadratic function in terms of the dimension
of a certain subspace of $\F_{p^n}$ (seen as a vector space over
$\F_p$) is essentially Proposition 2 in \cite{hk}. The result is
obtained via the standard squaring technique. We present the proof
to keep the paper self contained.
\begin{theorem}
\label{hk-prop} Let $f$
be the quadratic $p$-ary function
\[ f(x) = \Tr\left(\sum_{i=0}^la_ix^{p^i+1}\right), \]
and let $L(z)$ be the linearized polynomial
\begin{equation}
\label{L} L(z)=
\sum_{i=0}^l\left(a_i^{p^l}z^{p^{l+i}}+a_i^{p^{l-i}}z^{p^{l-i}}\right).
\end{equation}
The square of the Fourier transform of $f$ takes absolute values $0$
and $p^{n+s}$, where $s$ is the dimension of the kernel of the
linear transformation on $\F_{p^n}$ defined by $L(z)$.
\end{theorem}

{\it Proof.} With the standard squaring technique we obtain
\begin{eqnarray*}
|\widehat{f}(-b)|^2 & = &
\sum_{x,y\in\F_{p^n}}\epsilon_p^{f(x)-f(y)+\Tr(b(x-y))}\\ &=&
\sum_{y,z\in\F_{p^n}}\epsilon_p^{f(y+z)-f(y)+\Tr(bz)} \\
& = & \sum_{z\in\F_{p^n}}\epsilon_p^{f(z)+\Tr(bz)}
\sum_{y\in\F_{p^n}}\epsilon_p^{f(y+z)-f(y)-f(z)}.
\end{eqnarray*}
Observe that 
\begin{eqnarray*}
f(y+z) - f(y) - f(z)& = &
\Tr\left(\sum_{i=0}^la_i\left((y+z)^{p^i+1}-y^{p^i+1}-z^{p^i+1}\right)\right)\\
&=& \Tr\left(\sum_{i=0}^la_i\left(yz^{p^i}+y^{p^i}z\right)\right) \\
& = & \Tr\left(y^{p^l}\sum_{i=0}^l\left(a_i^{p^l}z^{p^{l+i}} +
a_i^{p^{l-i}}z^{p^{l-i}}\right)\right)\\  &=& \Tr(y^{p^l}L(z)).
\end{eqnarray*}
Consequently
\begin{eqnarray}
\label{supp} \nonumber |\widehat{f}(-b)|^2 & = &
\sum_{z\in\F_{p^n}}\epsilon_p^{f(z)+\Tr(bz)}
\sum_{y^{p^l}\in\F_{p^n}}\epsilon_p^{\Tr(yL(z))} \\ \nonumber
&=& p^n\sum_{z\in\F_{p^n} \atop L(z)=0}\epsilon_p^{f(z)+\Tr(bz)} \\ 
& = & \left\{\begin{array}{lr}
p^{n+s} & \mbox{if}\;f(z)+\Tr(bz) \equiv 0\;\mbox{on}\; ker(L) \\
0 & \mbox{otherwise}\end{array}\right.
\end{eqnarray}
where in the last step we used that $f(z) + \Tr(bz)$ is linear on
the kernel of $L$.
\hfill$\Box$\\[.5em]
\bigskip

\subsection{Monomials}

Our next goal is to find quadratic near-bent functions that can be
used to construct bent functions as described in Theorem \ref{thm1}.
We might hope for a monomial function, but unfortunately these do not exist
as we now prove.

\begin{theorem}
\label{yok} Quadratic monomial near-bent functions $f(x) =
\Tr(ax^{p^r+1})$, $a\in\F_{p^n}$, in odd characteristic $p$ do not
exist.
\end{theorem}
{\it Proof.} The linearized polynomial $(\ref{L})$ that corresponds
to $f(x) = \Tr(ax^{p^r+1})$ is given by $L(z) =
az+a^{p^r}z^{p^{2r}}$. We have to show that for any odd prime $p$,
integers $r,n \ge 1$ and $a \in \F_{p^n}$ the kernel of the linear
map on $\F_{p^n}$ induced by $L$ does not have dimension $1$. For a
primitive element $\gamma$ of $\F_{p^n}$ let $a = \gamma^c$ for some
$c, 0 \le c \le p^n-2$. Then $L(\gamma^t) = 0$ for an exponent $t,0
\le t \le p^n-2,$ if and only if
\[ \gamma^{\frac{p^n-1}{2}-c(p^r-1)} = \gamma^{(p^{2r}-1)t}, \]
which is equivalent to
\[ \frac{p^n-1}{2}-c(p^r-1) \equiv (p^{2r}-1)t \bmod (p^n-1). \]
Clearly the kernel of $L$ has dimension $1$ if and only if this
congruence has $p-1$ incongruent solutions, which applies if and
only if the two conditions $\gcd(p^{2r}-1,p^n-1) = p-1$ and $p-1$
divides $\frac{p^n-1}{2}-c(p^r-1)$ hold. The first condition is
satisfied if and only if $\gcd(2r,n) = 1$, in particular $n$ is then
odd, which contradicts the second condition.
\hfill$\Box$\\[.5em]

\begin{remark}
In \cite{hk} it is pointed out that $f$ is bent, i.e. the kernel of
$L$ has dimension $0$, if and only if $p^{\gcd(2r,n)}-1$ does not
divide $\frac{p^n-1}{2}-c(p^r-1)$.
\end{remark}

\subsection{Binomials}

As a consequence of Theorem \ref{yok} we must consider
non-monomial quadratic functions in order to be able to apply
Theorem \ref{thm1}. Two classes of binomial near-bent functions are
presented in the following theorem.
%
\begin{theorem}
\label{bino} Let $c\ne 0$ be an element of $\F_p$. 
The function $f$ from $\F_{p^n}$ to $\F_p$ given by
\begin{itemize}
\item[(i)] 
\begin{equation}
\label{binomial1} f(x) = \Tr\left(cx^{p^r+1}-cx^{p^t+1}\right)
\end{equation}
is near-bent if and only if $\gcd(n,r+t) = \gcd(n,r-t) = \gcd(n,p) = 1$,
and
\item[(ii)] 
\begin{equation}
\label{binomial} f(x) = \Tr\left(cx^{p^r+1}+cx^{p^t+1}\right)
\end{equation}
is near-bent if and only if $\gcd(n,2(r+t)) = \gcd(n,2(r-t)) = 2$, $r-t$ is odd,
and $\gcd(n,p) = 1$.
\end{itemize}
\end{theorem}
{\it Proof.} By Theorem \ref{hk-prop} a function $f$ from $\F_{p^n}$ to $\F_p$ is
near-bent if and only if the kernel (in $\F_{p^n}$) of the corresponding linearized 
polynomial $L(x)$ 
has dimension $1$ as a vector space over $\F_p$, i.e., $\gcd(L(x),x^{p^n}-x)$ has
degree $p$. Equivalently $ker(L)$ is one-dimensional if and only if the
associates $A(x)$ and $x^n-1$ of $L(x)$
%
%
and $x^{p^n}-x$, respectively, satisfy $\deg(\gcd(A(x),x^n-1)) = 1$, see
\cite[p.118]{ln}. \\
For the binomial $(\ref{binomial1})$ we have $L(x) = c(x+x^{p^{2r}}-x^{p^{r-t}}-x^{p^{r+t}})$,
consequently $A(x) = c(1+x^{2r}-x^{r-t}-x^{r+t}) = c(x^{r+t}-1)(x^{r-t}-1)$. Using
$\gcd(x^m-1,x^n-1) = x^{\gcd(m,n)}-1$ we easily see that $\deg(\gcd(A(x),x^n-1)) = 1$
if and only if $\gcd(n,r+t) = \gcd(n,r-t) = \gcd(n,p) = 1$. The last condition prevents
$1$ from being a multiple root of $x^n-1$. \\
The polynomial $A(x)$ for the binomial $(\ref{binomial})$ is given by 
$A(x) = c(1+x^{2r}+x^{r-t}+x^{r+t}) = c(x^{r+t}+1)(x^{r-t}+1)$. Using $\gcd(x^m-1,x^n-1) = 
x^{\gcd(m,n)}-1$ we obtain that
\begin{equation}
\label{gcd}
g=\gcd(x^{r\pm t}+1,x^n-1) = (x^{\gcd(2(r\pm t),n)}-1)/(x^{\gcd(r\pm t,n)}-1).
\end{equation}
If $n$ is odd then $g = 1$, thus we need $n$ even and hence we have $\gcd(n,2(r\pm t)) \ge 2$.
If $\gcd(n,2(r\pm t)) = 2$ then $r\pm t$ odd is a necessary and sufficient condition for
$g = x+1$. If $\gcd(n,2(r\pm t)) = 2u$ for an odd integer $u>1$ then by equation
$(\ref{gcd})$ we have $g = x^u+1$, and $\gcd(n,2(r\pm t)) = 2e$ for an even integer
$e$ yields $g = 1$ or $g = x^e+1$. 
As a consequence, $\gcd(n,2(r+t)) = \gcd(n,2(r-t)) = 2$ and $r-t$ ($r+t$) odd are necessary conditions
for $\gcd(A(x),x^n-1) = x+1$. As $-1$ is then a double root of $A(x)$ we
obtain $\gcd(A(x),x^n-1) = x+1$ if and only if $\gcd(n,p) = 1$.
\hfill$\Box$\\[.5em]
\begin{remark}
\label{rem1} The kernel of $L(x)$ in $\F_{p^n}$ for the function $(\ref{binomial1})$  is 
the set of the
solutions of $x^p-x$, which is $\F_p$.
For the function $(\ref{binomial})$ the kernel of $L(x)$ in $\F_{p^n}$ is the set of the
solutions of $x^p+x$, which are all in $\F_{p^2}$.
\end{remark}

\begin{remark}
Note that in part (i), $n$ can be either even or odd, whereas in part (ii),
the conditions imply that $n$ must be even.
\end{remark}

\section{Constructions of Bent Functions, Examples}
\label{sec4}

In order to apply Theorem \ref{thm1} we need $p$ near-bent functions
such that the supports of their Fourier transforms are pairwise
disjoint. We observe that the support of the Fourier transform of the
quadratic $p$-ary function $f(x) = \Tr(\sum_{i=0}^la_ix^{p^i+1})$ is
explicitly described in equation $(\ref{supp})$. This will be used
in the following theorem which describes how to obtain a set of
$p$ quadratic near-bent functions with the required properties.
%
\begin{theorem}
\label{prop3} Let $g_0,g_1,\ldots,g_{p-1}$ be quadratic near-bent
functions from $\F_{p^n}$ to $\F_p$ such that the linearized
polynomials $L_0,L_1,\ldots,L_{p-1}$ corresponding to
$g_0,g_1,\ldots,g_{p-1}$, respectively, have the same kernel
$\{c\beta, 0 \le c \le p-1\}$ in $\F_{p^n}$. Let
$b_0,b_1,\ldots,b_{p-1} \in \F_{p^n}$ such that
$g_k(\beta)+\Tr(b_k\beta) = g_0(\beta)+k$, $0 \le k \le p-1$. The
$p$ near-bent functions $f_0,f_1,\ldots,f_{p-1}$ defined by 
$f_k(x)= g_k(x)+\Tr(b_kx)$, $0 \le k \le p-1$, satisfy
$supp(\widehat{f_i})\cap supp(\widehat{f_j})=\emptyset$ for $0 \le
i\neq j \le p-1$.
\end{theorem}
{\it Proof.} We have to show that $-b \in supp(\widehat{f_j})$
implies $-b \not\in supp(\widehat{f_i})$ for integers $0 \le j,i \le
p-1$, $j \ne i$. Suppose $-b \in supp(\widehat{f_j})$, i.e.
$g_j(\beta) + \Tr(b_j\beta) + \Tr(b\beta) = g_0(\beta) + j +
\Tr(b\beta) = 0$.
Then $f_i(\beta) + \Tr(b\beta) = g_i(\beta) + \Tr(b_i\beta) +
\Tr(b\beta) =
g_0(\beta) + i + \Tr(b\beta) \ne 0$ when $i \ne j$. \hfill$\Box$\\[.5em]
{\it Example 1.} Let $g(x) = \Tr(\sum_{i=0}^la_ix^{p^i+1})$ be a
quadratic near-bent function from $\F_{p^n}$ to $\F_p$, let
$\beta\in \F_{p^n}^*$ be a root of the corresponding linearized
polynomial $L(x)$ and let $b_k \in \F_{p^n}$, $0 \le k \le p-1$,
such that $\Tr(b_k\beta) = k$. Then the function $F_1(x,y)$ from
$\F_{p^n}\times\F_p$ to $\F_p$ given by
\[ F_1(x,y)=(p-1)\sum_{k=0}^{p-1}\frac{y(y-1)\cdots(y-(p-1))}{y-k}
\left(\Tr(\sum_{i=0}^la_ix^{p^i+1})+\Tr(b_kx)\right) \] is bent. 
Here all the $g_k(x)$ are equal to $g(x)$.
As
easily observed the bent function $F_1(x,y)$ is quadratic, a result
of the fact that the $p$ near-bent functions 
\[
f_k(x)=\Tr(\sum_{i=0}^la_ix^{p^i+1})+\Tr(b_kx)
\]
used for the
construction only differ in a linear term. 

As the values of the
Fourier spectrum of this quadratic bent function are the nonzero
values in the Fourier spectrum of the underlying near-bent function,
$(\ref{(2)})$ implies that the nonzero normalized Fourier coefficients
of a quadratic near-bent function $f$ from $\F_{p^n}$ to $\F_p$
satisfy
\[ p^{-n/2}\widehat{f}(b) =
\left\{\begin{array}{r@{\quad:\quad}l}
\pm \epsilon_p^{J(b)} &  n\;\mbox{odd and}\;p \equiv 3\bmod 4 \\
\pm i\epsilon_p^{J(b)} & n\;\mbox{even or}\;n\;\mbox{odd and}\;p
\equiv 1\bmod 4
\end{array}\right. \]
for some function $J(x)$ from $supp(\widehat{f})$ to $\F_p$. 
By \cite[Proposition 1]{hk} quadratic bent functions are always
(weakly) regular, and thus  quadratic near-bent functions are also, in
that sense that $\zeta p^{-n/2}\widehat{f}(b) = \epsilon_p^{J(b)}$ for
all $b \in supp(\widehat{f})$ and a fixed complex number $\zeta$ with
absolute value $1$ (in this connection we remark that adding a
linear term to a $p$-ary function $f$ does not change the Fourier
spectrum). A detailed description of the Fourier spectrum of quadratic
near-bent functions will be given in Theorem \ref{Fourier}.\\

As by Remark \ref{rem1} the linearized polynomials of all near-bent
functions of the form $x^{p^r+1}+x^{p^t+1}$ have the same kernel 
(the solutions of $x^p+x$) we now apply Theorem
\ref{prop3} to construct non-quadratic bent functions. \\

{\it Example 2.} 
Let $p=3$ and $n=8$.  We construct a bent function in dimension 9.
Let $g_0(x)=g_1(x)=\Tr(x^{3^2+1}+x^{3+1}), g_2(x)=\Tr(x^{3^6+1}+x^{3^5+1})$ 
be functions from $\F_{3^8}$ to $\F_3$. By the remark after Theorem \ref{bino}
the kernel $ker(L)$ of the corresponding linear transformations
consists of the solutions of $x^3+x$. For a root $\beta$ of $x^2+1$
we have $g_0(\beta)=g_1(\beta)=g_2(\beta)=0$ and therefore
we need to find $b_0,b_1,b_2\in \F_{3^8}$ such that ${\rm
Tr_8}(b_j\beta)=j$, to construct three near-bent functions
such that the supports of their Fourier transforms are pairwise
disjoint. Observing that ${\rm Tr_8}(\beta)=0, {\rm
Tr_8}(\beta^2)=1$ and ${\rm Tr_8}(2\beta^2)=2$, we can choose
$b_0=1,b_1=\beta, b_2=2\beta$, and we therefore set $f_0(x)={\rm
Tr_8}(x^{10}+x^4+b_0x), f_1(x)={\rm Tr_8}(x^{10}+x^4+b_1x),
f_2(x)={\rm Tr_8}(x^{3^6+1}+x^{3^5+1}+b_2 x)$. By Theorem 1, the
following function from $\F_{3^8}\times\F_3$ to $\F_3$ of algebraic
degree $4$ is bent:
\begin{eqnarray*}
&  & 2\sum_{k=0}^2\frac{y(y-1)(y-2)}{y-k}f_k(x) \\
& = & (2y^2+1){\rm Tr_8}(x^{10}+x^4+x)+(2y^2+2y){\rm
Tr_8}(x^{10}+x^4+\beta x)+ \\
& & (2y^2+y){\rm Tr_8}(x^{3^6+1}+x^{3^5+1}+2\beta x) \\
& = & 2y^2{\rm Tr_8}(2x^{10}+2x^4+x^{3^6+1}+x^{3^5+1}+x)+ \\
& & y{\rm Tr_8}(2x^{10}+2x^4+x^{3^6+1}+x^{3^5+1}+\beta x)+{\rm
Tr_8}(x^{10}+x^4+x) \\
& = & 2y^2{\rm Tr_8}(2x^4+x^{3^5+1}+x)+ y{\rm
Tr_8}(2x^4+x^{3^5+1}+\beta x) \\
& & +{\rm Tr_8}(x^{10}+x^4+x).
\end{eqnarray*}

\vspace{.5em}
{\it Example 3.} 
Let $g_0(x)$, $g_1(x)$, $g_2(x)$, $b_0, b_1, b_2$ be as
in Example 2, except use $2g_1(x)$ in place of $g_1(x)$.

\vspace{.5em}
{\it Example 4.} 
Let $g_0(x)$, $g_1(x)$, $g_2(x)$, $b_0, b_1, b_2$ be as
in Example 2, except take $g_1(x)=g_2(x)$ instead of taking $g_1(x)$
to be $g_0(x)$.

\vspace{.5em}
{\it Example 5.} 
Let $g_0(x)$, $g_1(x)$, $g_2(x)$, $b_0, b_1, b_2$ be as
in Example 4, except use $2g_1(x)$ in place of $g_1(x)$.

\vspace{.5em}
{\it Example 6.} 
We give an example in even dimensions using
Theorem \ref{bino} part (i); let $p=3$ and $n=5$.
Let $g_0(x)=\Tr(x^{3^2+1}-x^{3+1}), g_1(x)=\Tr(2x^{3^2+1}-2x^{3+1})$,
$g_2(x)=\Tr(x^{3^2+1}-x^{3+1})$ be functions from $\F_{3^5}$ to $\F_3$. 
Then these functions vanish on $\F_3$, so 
we can choose $b_0=0,b_1=2, b_2=1$ which yields
$f_0(x)={\rm Tr_5}(x^{10}-x^4)$, $f_1(x)={\rm Tr_5}(2x^{10}-2x^{4}+2 x)$,
$f_2(x)={\rm Tr_5}(x^{10}-x^{4}+ x)$. 
The resulting bent function 
\[
 2\sum_{k=0}^2\frac{y(y-1)(y-2)}{y-k}f_k(x) 
 =2(y-1)(y-2)f_0(x)+2y(y-2)f_1(x)+2y(y-2)f_2(x)
\]
from 
$\F_{3^5}\times\F_3$ to $\F_3$ again has degree 4.

\begin{remark}
We will see later that Examples 2, 4, 6 are weakly regular, but examples
3 and 5 are not weakly regular.
\end{remark}

\section{(Non) Weak Regularity}\label{sec5}

We finally consider the question of whether the bent functions obtained
from quadratic near-bent functions with Theorem \ref{thm1}, Theorem
\ref{bino} and Theorem \ref{prop3} are weakly regular. 
In this section we will prove necessary and sufficient conditions for weak regularity.
Throughout this section, $\eta$ denotes the quadratic
character in $\F_p$.

\subsection{Necessary and sufficient
conditions for weak regularity, and the Fourier spectrum}\label{necsuf}

We start with explicitly determining the Fourier spectrum of a quadratic
near-bent function $f$. Choosing and fixing a basis
$\{\alpha_1,\ldots,\alpha_n\}$ of $\F_{p^n}$ over $\F_p$ we
correspond $x = \sum_{i=1}^nx_i\alpha_i$ to the vector $\mathbf{x} =
(x_1,\ldots,x_n)$. Then we can associate a quadratic function $f(x)
= \Tr(\sum_{i=0}^la_ix^{p^i+1})$ with a quadratic form
\[ f(\mathbf{x}) = \mathbf{x}^TA\mathbf{x} \]
where $\mathbf{x}^T$ denotes the transpose to the vector
$\mathbf{x}$, and
the matrix $A$ has entries in $\F_p$.
By \cite[Theorem 6.21]{ln} any quadratic form is
equivalent to a diagonal quadratic form, i.e. $D=C^TAC$ for a
nonsingular matrix $C$ over $\F_p$ and a diagonal matrix $D =
diag(d_1,\ldots,d_n)$. With standard arguments based on Theorems
6.26 and 6.27 in \cite{ln} one can express the Fourier transform of a
quadratic near-bent function in terms of the product of the nonzero
entries in $D$ (for bent functions see \cite[Proposition 1]{hk}).
\begin{theorem}
\label{Fourier} Let $f$ be a quadratic near-bent function from
$\F_{p^n}$ to $\F_p$ and $f(\mathbf{x}) = \mathbf{x}^TA\mathbf{x}$ be
the associated quadratic form. Then a corresponding diagonal matrix
$D$ has $n-1$ (not necessarily distinct) nonzero entries
$d_1,\ldots,d_{n-1}$, and the Fourier spectrum of $f$ is given by
\[
\begin{array}{r@{\quad:\quad}l}
\left\{0, \eta(\Delta)p^{\frac{n+1}{2}}\epsilon_p^{J(b)}\right\} & p\equiv 1\bmod 4, \\[.5em]
\left\{0,
(-1)^{\frac{n-2}{2}}\eta(\Delta)ip^{\frac{n+1}{2}}\epsilon_p^{J(b)}\right\}
&
p\equiv 3\bmod 4\;\mbox{and}\;n\;\mbox{even},\\[.5em]
\left\{0,
(-1)^{\frac{n-1}{2}}\eta(\Delta)p^{\frac{n+1}{2}}\epsilon_p^{J(b)}\right\}
& p\equiv 3\bmod 4\;\mbox{and}\;n\;\mbox{odd},
\end{array}
\]
where $J(x)$ is a function from $supp(\widehat{f})$ to $\F_p$,
$\Delta = \prod_{i=1}^{n-1}d_i$, and $\eta$ denotes the quadratic
character in $\F_p$.
\end{theorem}
{\it Proof.} We write $f(x)-\Tr(bx) = j$ equivalently as
\begin{equation}
\label{=j} \mathbf{x}^TA\mathbf{x} -
(\Tr(b\alpha_1),\ldots,\Tr(b\alpha_n))\mathbf{x} =
\mathbf{x}^TA\mathbf{x} + c^T\mathbf{x} = j
\end{equation}
where $c \in \F_p^n$. Denoting by $N_b(j)$, $j=0,\ldots,p-1$, the
number of solutions in $\F_p^n$ for $(\ref{=j})$ we observe that
\[ \widehat{f}(b) = \sum_{j=0}^{p-1}N_b(j)\epsilon_p^j. \]
Substituting $\mathbf{x} = C\mathbf{y}$ where $C$ is the nonsingular
matrix with $D = C^TAC$ we obtain for equation $(\ref{=j})$
\[ \mathbf{y}^TD\mathbf{y} + c^TC\mathbf{y} = \sum_{i=1}^n(d_iy_i^2+c_iy_i) = j \]
with $c^TC = (c_1,\ldots,c_n)^T \in \F_p^n$. Suppose that w.l.o.g.
$d_n = 0$ is the only zero in the diagonal of $D$. Performing the
substitution $y_i = z_i-c_i/(2d_i)$ for $i = 1,\ldots,n-1$, we
get
\begin{equation}
\label{=ja} \sum_{i=1}^{n-1}d_iz_i^2 = j +
\sum_{i=1}^{n-1}\frac{c_i^2}{4d_i} - c_ny_n.
\end{equation}
We note that finding solutions $(z_1,\ldots,z_{n-1},y_n) \in \F_p^n$
for $(\ref{=ja})$ is equivalent to finding solutions
$(x_1,\ldots,x_n) \in \F_p^n$ for $(\ref{=j})$,
and the number of solutions is $N_b(j)$.
Also note that the map $b\mapsto (c_1,\ldots,c_n)$ is a bijection.

First suppose we have  $b$ with $c_n \ne 0$. 
For an arbitrary choice of
$j,z_1,\ldots,z_{n-1}$, equation $(\ref{=ja})$ is satisfied for a
unique choice for $y_n$. As a consequence
$N_b(j)$ has the same value for each $j$  and so $\widehat{f}(b) = 0$. 
(As an aside, there are $p^{n-1}(p-1)$ vectors $(c_1,\ldots,c_n)$ with $c_n\ne 0$, so
$p^{n-1}(p-1)$ is the multiplicity of 0 in the Fourier spectrum.)

Now suppose we have  $b$ such that $c_n = 0$, and define $J(b) =
-\sum_{i=1}^{n-1}\frac{c_i^2}{4d_i}$.
We need to consider the cases of even and odd $n$ separately.
For even $n$, \cite[Theorem 6.27]{ln} gives the number of
solutions of (\ref{=ja}) (i.e. gives $N_b(j)$) and we have
\begin{eqnarray}
\label{evenn} \nonumber \widehat{f}(b) & = &
p\sum_{j=0}^{p-1}\left(p^{n-2}+
p^{(n-2)/2}\eta\left((-1)^{(n-2)/2}(j-J(b))\Delta\right)\right)\epsilon_p^j \\
& = &
p^{n/2}\eta\left((-1)^{(n-2)/2}\Delta\right)\sum_{j=0}^{p-1}\eta(j-J(b))\epsilon_p^j.
\end{eqnarray}
By \cite[Theorem 5.15]{ln} we have
\[ \sum_{j=0}^{p-1}\eta(j)\epsilon_p^j = \left\{
\begin{array}{r@{\quad:\quad}l}
p^{1/2} & p \equiv 1 \bmod 4, \\
ip^{1/2} & p \equiv 3 \bmod 4,
\end{array}\right. \]
(with the usual convention that $\eta(0) = 0$), thus equation
$(\ref{evenn})$ reduces to
\[ \widehat{f}(b) = (-1)^{\frac{(p-1)(n-2)}{4}}\eta(\Delta)i^{s(p)}p^{(n+1)/2}\epsilon_p^{J(b)}, \]
where $s(p) = 0$ if $p \equiv 1 \bmod 4$ and else $s(p) = 1$. 

For odd $n$, Theorem 6.26 in \cite{ln} implies that the values
$N_b(j)$ are all equal except for $N_b(J(b))$ which differs from
the others by $\eta((-1)^{(n-1)/2}\Delta)p^{(n+1)/2}$. Consequently
\[ \widehat{f}(b) = (-1)^{\frac{(p-1)(n-1)}{4}}\eta(\Delta)p^{(n+1)/2}\epsilon_p^{J(b)}, \]
which shows the correctness of the values for the Fourier transform
given in the theorem.

Finally, as shown in \cite{hk}, $f$ is bent if and only if the
associated quadratic form is nondegenerate. If the rank of $A$ is
$n-s$, and w.l.o.g. $d_1,\ldots,d_{n-s} \ne 0$, we consider elements
$b \in \F_{p^n}$ for which $c_{n-s+1}=\cdots =c_n = 0$
(if one is nonzero, then $\widehat{f}(b) = 0$).
Here we use
that $C$ is nonsingular and that $b \mapsto
(\Tr(b\alpha_1),\ldots,\Tr(b\alpha_n))$ defines a one-to-one linear
transformation from $\F_{p^n}$ to $\F_p^n$. Then with the same
arguments as above we obtain that
$|\widehat{f}(b)| = p^{(n+s)/2}$. Thus $f$ is not near-bent for $s \ne 1$. \hfill$\Box$\\

\begin{remark}
We remark that the `aside' comment in the previous proof gives
a description of the $b$ with $\widehat{f}(b) = 0$.
\end{remark}

As an immediate consequence we obtain the following corollary, which
gives our necessary and sufficient conditions for weak regularity.
\begin{corollary}
\label{cor2} Let $f_0,\ldots,f_{p-1}$ be $p$-ary quadratic near-bent
functions with $supp(\widehat{f_i})\cap supp(\widehat{f_j})
=\emptyset$ for $0\le i\ne j \le p-1$. Let $A_i$, $0 \le i \le p-1$,
be the matrix of the quadratic form associated with $f_i$, and let
$\Delta_i$ be the product of the nonzero eigenvalues of $A_i$,
respectively. Then the bent function constructed as in Theorem
$\ref{thm1}$ is weakly regular if and only if $\eta(\Delta_0) =
\eta(\Delta_1) = \cdots = \eta(\Delta_{p-1})$.
\end{corollary}

\subsection{(Non) Weak regularity of our examples}

In order to decide whether the bent functions obtained 
in Theorem \ref{prop3} using the functions of Theorem \ref{bino}
are weakly regular, we are interested in the matrices (and their
eigenvalues) of the quadratic forms associated with these quadratic
functions. We start with a lemma on the quadratic character of the
product of the nonzero eigenvalues of the matrix $A$ associated with
a quadratic near-bent function.
\begin{lemma}
\label{sign}
Let $f$ be a quadratic near-bent function from $\F_p^n$ to $\F_p$, let $\mathbf{x}^TA\mathbf{x}$ 
be the associated quadratic form and let $\Delta$ denote the product of the $n-1$
nonzero eigenvalues of $A$. For a nonzero constant $c \in \F_p$, the product $\Delta^{(c)}$
of the $n-1$ nonzero eigenvalues of the matrix
for the quadratic form associated with $cf$ satisfies
\[ \eta(\Delta^{(c)}) = \eta(c)^{n-1}\eta(\Delta). \]
\end{lemma}

\bproof 
The lemma follows with $cf(\mathbf{x}) = \mathbf{x}^TcA\mathbf{x}$.
\eproof \\[.5em]

We observe that as a consequence of Theorem \ref{Fourier} and Lemma \ref{sign},
if $n$ is even then we can change the signs of the Fourier coefficients if we switch
from the near-bent function $f$ to the near-bent function $cf$, with a nonsquare
$c \in \F_p$. With these observations we can obtain an infinite class of non-weakly
regular bent functions. As building blocks we may use the near-bent functions from
Theorem \ref{prop3}.

\begin{theorem}
\label{non-weakly}
Let $n$ be even.
For each $0\le k \le p-1$ let $c_k$ be an element of $\F_p^*$ and let $g_k$ be the near-bent function
\[ g_k(x) = c_k\Tr(x^{p^{r}+1}-x^{p^{t}+1})\quad\mbox{or}\quad g_k(x) = c_k\Tr(x^{p^{r}+1}+x^{p^{t}+1}) \]
from $\F_{p^n}$ to $\F_p$ described as in Theorem $\ref{bino}(i)$ or Theorem $\ref{bino}(ii)$,
respectively. Here $r=r_k$ and $t=t_k$ may vary with $k$. 
Let $f_0,\ldots,f_{p-1}$ be $p$-ary quadratic near-bent functions
from $\F_{p^n}$ to $\F_p$ with $supp(\widehat{f_i})\cap
supp(\widehat{f_j}) =\emptyset$ for $i\ne j$, which are
obtained as in Theorem $\ref{prop3}$ from the binomial
near-bent functions $g_k$. Then the bent function $F(x,y)$ 
from $\F_{p^n}\times\F_p$ to $\F_p$
\[ F(x,y)=(p-1)\sum_{k=0}^{p-1}\frac{y(y-1)\cdots(y-(p-1))}{y-k}f_k(x) \]
is weakly regular for $(p-1)^p/2^{p-1}$ choices for
$(c_0,\ldots,c_{p-1}) \in (\F_p^*)^{p}$, and non-weakly regular for the
remaining $(2^{p-1}-1)(p-1)^p/2^{p-1}$ choices for $(c_0,\ldots,c_{p-1}) \in (\F_p^*)^{p}$.
\end{theorem}

\bproof 
Let $\Delta_k$ denote the product of the nonzero eigenvalues of the matrix that corresponds
to $g_k$. Then by Corollary \ref{cor2} and Lemma \ref{sign}, $F(x,y)$ is weakly regular if 
and only if $\eta(c_0)^{n-1}\eta(\Delta_0) = \eta(c_1)^{n-1}\eta(\Delta_1)
= \cdots = \eta(c_{p-1})^{n-1}\eta(\Delta_{p-1})$.
Since $n-1$ is odd, for a fixed value of $\eta(c_k)^{n-1}\eta(\Delta_k)$, $1$ or $-1$,
for every $0 \le k \le p-1$ we have $(p-1)/2$ choices for $c_k$.
This gives in total $2((p-1)/2)^p$ choices for $(c_0,\ldots,c_{p-1}) 
\in (\F_p^*)^{p}$ for which $F(x,y)$ is weakly regular.
\eproof \\[.5em]

The conditions in the previous theorem mean that a simple tweak of the coefficients 
$c_k$ can drastically change the nature of the bent function. For example,
we observe that the bent functions from Example 2 and 4 where each $c_k$ is equal to 1,
are weakly regular. The normalized Fourier spectrum (without multiplicities) of 
Example 2 (and 4) is shown in Fig 4, and the spectrum with multiplicities is
\[
( -i)^{2187}, (-i \epsilon_3 )^{2268},  (-i \epsilon_3^2)^{2106}, 
\]
(recall the dimension is 9).
In Examples 3 and 5, still $c_0 = c_2 = 1$, but we changed $c_1$ to $c_1 = 2$.
This alters the sign of the Fourier coefficients for one of the near-bent functions.   
Consequently the bent functions in Examples 3 and 5 are non-weakly regular.
The normalized Fourier spectrum (without multiplicities) of Example 3 (and 5) is 
shown in Fig 6, and
the spectrum with multiplicities is
\[
(i \epsilon_3 )^{702}, (i \epsilon_3^2)^{756}, ( -i)^{1458}, i^{729}, (-i \epsilon_3^2)^{1404}, 
(-i \epsilon_3 )^{1512}.
\]
Many non-weakly regular bent functions can clearly be constructed in this manner;
one simply has to arrange that for two of the coefficients $c_k$ and $c_{k^\prime}$
we have $\eta(c_k)\eta(\Delta_k) \ne \eta(c_{k^\prime})\eta(\Delta_{k^\prime})$.

Finally we remark that the arguments in Theorem \ref{non-weakly} are not restricted to 
near-bent functions of the form $(\ref{binomial1})$, $(\ref{binomial})$, but applicable 
to every set of $p$ quadratic near-bent functions in even dimension satisfying the 
conditions of Theorem \ref{thm1}.

\subsection{A Family of Weakly Regular Bent Functions}

We show that the bent functions obtained with Theorem \ref{thm1} using the
near-bent functions $(\ref{binomial1})$ are always weakly regular when $n$ is odd.
As weakly regular bent functions are useful for the construction of certain combinatorial
objects such as partial difference sets, strongly regular graphs and association schemes 
(see \cite{tpf}) they are of independent interest. 

First observe that a quadratic function $f(x) = \Tr(\sum_{i=0}^la_ix^{p^i+1})$ can be written  
in the form $f(x) = \Tr(x\mathcal{L}(x))$ for the linear transformation
$\mathcal{L}(x) = \sum_{i=0}^la_ix^{p^i}$ on $\F_{p^n}$. 
As we did in Section \ref{necsuf}, if
we choose a basis $\{\alpha_0,\ldots,\alpha_{n-1}\}$ of $\F_{p^n}$ over $\F_p$, 
and let $x = \sum_{i=0}^{n-1}x_i\alpha_i \in \F_{p^n}$ correspond 
to the vector $\mathbf{x} = (x_0,\ldots,x_{n-1})$, then
we can associate $f(x)$ with a quadratic form $\mathbf{x}^TA\mathbf{x}$.

We now choose $\{\alpha_0,\ldots,\alpha_{n-1}\}$ to
be a self-dual basis of $\F_{p^n}$ over $\F_p$, which exists if and only if $n$ is odd.
Then  a straightforward calculation shows that $A$ is the
matrix representation of the linear transformation induced by $\mathcal{L}(x)$ with respect to
the given self-dual basis.

For the functions of interest to us, $\mathcal{L}(x)$ is of
the form
\[ \mathcal{L}(x) = cx^{p^r}-cx^{p^t} \]
where $c\in \F_p$. We choose the self-dual basis to be also a normal basis
$\{\alpha,\alpha^p,\ldots,\alpha^{p^{n-1}}\}$ (which always is possible when $n$ is odd,
see \cite{ln}).
Then $\mathcal{L}(\alpha_i) = \mathcal{L}(\alpha^{p^i}) =
c\alpha^{p^{r+i}}-c\alpha^{p^{t+i}} = c\alpha_{r+i\bmod
n}-c\alpha_{t+i\bmod n}$. 
Hence the corresponding matrix $A_c^{(r,t)}$
is an $n\times n$ circulant matrix with first column
$(a_0,\ldots,a_{n-1})^T$ with $a_r = c, a_t = -c$ and $a_i = 0, i \ne
r,t$, or equivalently with first row $(s_0,\ldots,s_{n-1})$ with
$s_i=c$ if $i \equiv n-r \bmod n$, $s_i = -c$ if $i \equiv n-t \bmod n$, 
and $s_i=0$ otherwise. Using the notation
$A=C(s_0,\ldots,s_{n-1})$ for a circulant matrix $A$ with first row
$(s_0,\ldots,s_{n-1})$ we can summarize these observations with
\begin{eqnarray}
\label{matrix} \nonumber
A_c^{(r,t)} & = & C(s_0,\ldots,s_{n-1})\quad\mbox{where} \\
s_i & = & \left\{
\begin{array}{r@{\quad:\quad}l}
c & i \equiv n-r \bmod n \\
-c & i \equiv n-t \bmod n \\
0 & \mbox{otherwise.}
\end{array}\right.
\end{eqnarray}
We will use the following result on eigenvalues of circulant
matrices, see 1.6 in \cite{w}.
\begin{lemma}
\label{circeig} Let $n$ be an integer relatively prime to $p$, $u$ a
primitive $n$-th root of unity over $\F_p$, and let
$A=C(s_0,\ldots,s_{n-1})$ be an $n\times n$ circulant matrix. The
eigenvalues of $A$ are given by
\[ \lambda_j = \sum_{i=0}^{n-1}s_iu^{ij},\; j = 0,\ldots, n-1. \]
\end{lemma}
We denote the eigenvalues of the matrix $A_1^{(r,t)}$ corresponding to
$x^{p^r}-x^{p^t} \in \F_{p^n}[x]$ 
as $\lambda_j^{(r,t)}$, $j=0,1,\cdots,n-1$. By $(\ref{matrix})$ and 
Lemma \ref{circeig} we then have $\lambda_j^{(r,t)}=u^{(n-r)j}-u^{(n-t)j}$
for $j=0,1, \cdots,n-1$, where $u$ is a primitive $n$-th root of
unity over $\F_p$. As easily seen $\lambda_0^{(r,t)} = 0$ and
hence the product of the nonzero eigenvalues of $A_1^{(r,t)}$ is
$\Delta^{(r,t)}=\prod_{j=1}^{n-1}\lambda_j^{(r,t)}$.
%
%
Next we show that $\Delta^{(r,t)}$ does not depend on the special choice of $r$ and $t$.
\begin{lemma}
\label{lemma}
For an odd integer $n$ let $r,t$ and $v,w$ be pairs of
integers satisfying the conditions in Theorem \ref{bino}(i). Then with
the above notations
\[ \Delta^{(r,t)} = \prod_{j=1}^{n-1}\lambda_j^{(r,t)} = 
\prod_{j=1}^{n-1}\lambda_j^{(v,w)} = \Delta^{(v,w)}. \]
\end{lemma}
\bproof For fixed $j,\; 0 \le j \le n-1$, we are interested in two
integers $0 \le k_j,c_j \le n-1$ such that
\[ \lambda_j^{(r,t)} = u^{(n-r)j}-u^{(n-t)j} = (u^{(n-v)k_j}-u^{(n-w)k_j})u^{c_j} = \lambda_{k_j}^{(v,w)}u^{c_j}. \]
We therefore consider the linear system in the two variables $k_j,
c_j$
\begin{eqnarray*}
-rj & \equiv & -vk_j+c_j \bmod n \\
-tj & \equiv & -wk_j+c_j \bmod n
\end{eqnarray*}
yielding
\[ k_j \equiv \frac{r-t}{v-w}j \bmod n \quad\mbox{and}\quad c_j \equiv -\frac{rw-tv}{v-w}j \bmod n. \]
We remark that $k_j$ and $c_j$ are well defined since $\gcd(v-w,n) =
1$. Moreover also $\gcd(r-t,n) = 1$ thus $(r-t)/(v-w)$ is an
invertible residue modulo $n$ and $k_j$ runs through the integers
modulo $n$ if $j$ does. Consequently
\begin{eqnarray*}
\Delta^{(r,t)} & = & \prod_{j=1}^{n-1}\lambda_j^{(r,t)}=
\prod_{j=1}^{n-1}(u^{(n-v)k_j}-u^{(n-w)k_j})u^{c_j} \\
& = & \prod_{j=1}^{n-1}\lambda_{k_j}^{(v,w)}\prod_{j=1}^{n-1}u^{c_j} 
= \Delta^{(v,w)}\prod_{j=1}^{n-1}u^{c_j}.
\end{eqnarray*}
Since
\[ \prod_{j=1}^{n-1}u^{c_j} = u^{\sigma\sum_{j=1}^{n-1}j} = 1, \]
where $\sigma \equiv -\frac{rw-tv}{v-w} \bmod n$, the proof is complete.
\eproof \\[.5em]
\begin{theorem}\label{sum}
Let $n$ be odd.
For each $0\le k \le p-1$ let $c_k$ be an element of $\F_p^*$ and let $g_k$ be the near-bent function
\[ g_k(x) = c_k\Tr(x^{p^{r}+1}-x^{p^{t}+1}) \]
from $\F_{p^n}$ to $\F_p$ described as in Theorem $\ref{bino}(i)$.
Here $r=r_k$ and $t=t_k$ may vary with $k$. 
Let $f_0,\ldots,f_{p-1}$ be $p$-ary quadratic near-bent functions
from $\F_{p^n}$ to $\F_p$ with $supp(\widehat{f_i})\cap
supp(\widehat{f_j}) =\emptyset$ for $i\ne j$, which are
obtained as in Theorem $\ref{prop3}$ from the binomial
near-bent functions $g_k$.  Then the bent function $F(x,y)$ 
from $\F_{p^n}\times\F_p$ to $\F_p$
\[ F(x,y)=(p-1)\sum_{k=0}^{p-1}\frac{y(y-1)\cdots(y-(p-1))}{y-k}f_k(x) \]
is weakly regular.
\end{theorem}
{\it Proof.}
By Lemma \ref{sign} and Lemma \ref{lemma} we may write $c_k^{n-1}\Delta$ for the 
product of the nonzero eigenvalues of the circulant matrix we correspond to $g_k$, $0 \le k \le p-1$. 
Since $n-1$ is even we have $\eta(c_k^{n-1}\Delta) = \eta(\Delta)$ for all $0 \le k \le p-1$.
By Corollary \ref{cor2} we obtain then the assertion of the theorem. 
\hfill$\Box$

\begin{remark}
Theorem \ref{sum} implies that Example 6 is a weakly regular bent function.
\end{remark}

%
%

\end{document}